\documentclass[]{amsart}

\theoremstyle{definition}

\theoremstyle{remark}

\numberwithin{equation}{section}

\begin{document}

\title{Clarkson-Erd\"{o}s-Schwartz  Theorem on  a  Sector}

\author{Guan-Tie DENG}
\address{Sch. Math. Sci. \& Lab. Math. Com.
Sys., Beijing  Normal University, Beijing 100875, China}
\email{denggt@bnu.edu.cn}
\thanks{This work was partially supported by
NSFC (Grant  11071020) and by SRFDP (Grant 20100003110004). }



\date{\today }



\begin{abstract}
 We prove a
Clarkson-Erd\"{o}s-Schwartz type theorem for the case of a closed
sector in the plane. Concretely, we get some sufficient conditions
for the incompleteness and minimality of a M\"{u}ntz system
$E(\Lambda)=\{z^{\lambda_n}:n=0,1,\cdots\}$ in the space
$H_\alpha$, where $H_\alpha=A(I_\alpha)$,
$I_\alpha=\{z\in\mathbb{C}:|\arg (z)|\leq \alpha \text{ and }
|z|\leq 1\}$ and $A(K)=C(K)\cap H(\textbf{Int}[K])$ denotes the
space of continuous functions on the compact set $K$ which are
analytic in the interior of $K$. Furthermore, we prove that, if
$\textbf{span}[E(\Lambda)]$ is not dense in $H_\alpha$ then all
functions $f\in \overline{\textbf{span}}[E(\Lambda)]$ can be
analytically extended to the interior of the sector $I_\pi$.
\end{abstract}

\maketitle

\section{ Introduction }
  Suppose  $ \Lambda =\{ \lambda_n:n=1,2, \cdots \} $ is  a
sequence of  positive real numbers arranged  for convenience in
non-decreasing order:
$$
 0=\lambda_0< \lambda_1\leq \lambda_2\leq \cdots
, \ \mbox{ and  }  \delta (\Lambda )= \inf \{ \lambda
_{n+1}-\lambda_{n}:n\geq 0 \}.
$$
$ E(\Lambda )=\{z^{\lambda _n} =\exp \{ \lambda _n \log
z\}:n=0,1,2, \cdots \}$ is called a M\"{u}ntz system, $
\textbf{span}[E(\Lambda)]$ denotes the linear span of the
M\"{u}ntz system. The elements of the set
$\textbf{span}[E(\Lambda)]$ are called the M\"{u}ntz polynomial or
the $\Lambda$-polynomials [4]. Let $K$ be a compact set in the
plane, and let $C(K)$ be the space of all continuous functions on
$K$, equipped with the uniform norm. The famous M\"{u}ntz theorem
([4] and [18]) states that
  $\textbf{span}[E(\Lambda )]$  is
dense in $C([0,1])$  if and only if $ \
 \ \sum_{n=1}^{\infty}1/\lambda_n\ $  diverges.  Moreover, If the
set $\textbf{span}[E(\Lambda)]$ is not a dense subspace of
$C[0,1]$, it is natural to ask for a characterization of the
elements of its topological closure. This problem was  solved by
Clarkson and Erd\"{o}s [5] for the case of
 integer exponents $ \Lambda $, they proved that if  $ \
 \ \sum_{n=1}^{\infty}1/\lambda_n\ $  converges,  the elements in  the closure
   of  $\textbf{span}[E(\Lambda)]$
      can be extended
analytically throughout to the unit disc  with a series expansion
of the form $$
 f(x)= \sum_{k=0}^{\infty}a_kx^{\lambda_k}, \ 0\leq x<1  .
\eqno(1)
$$
This same question was also tackled by L.Schwarts [16] for certain
strictly increasing sequences of exponents (he assumed $ \delta
(\Lambda )>0$)  and by Borwein [4] and \'Erdelyi [6]. Nowadays
these results are referenced under the join name of
Clarkson-Erd\"{o}s-Schwartz Theorem [1].
  On the other hand, if $K$
is a compact in the plane whose complement is connected,
Mergelyan's theorem [15] claims  that the space  of complex
polynomials is a dense set of  $C(K)$. It is a nontrivial problem
to establish  a Clarkson-Erd\"{o}s-Schwartz  Theorem on a closed
sector in the plane. The aim of this paper  answer to this
problem. First we introduce some notations and definitions. Let
$B$ be a Banach space. If $E=\{ e_k:k=1,2,\cdots \}\subset B$, let
$\textbf{span}[E]$ denote the subspace of $B$, consisting of all
finite linear combinations of $E$ and let
$\overline{\textbf{span}}[E]$ be the closure of $\textbf{span}[E]$
in $B$. The set $ E $ is said to be {\it incomplete } in $B$ [17]
if $\overline{\textbf{span}}[E]$ does not coincide with the whole
$B $. The set $E $ is called to be ${\it minimal }$ in $B$ [17] if
no element of $E$ belongs to the closure of the vector subspace
generated by the other elements of $E$, i.e., for all $ e\in E $,
$e \not\in \overline{\textbf{span}}[E- \{e \}] $. The minimality
of the
 $E $ is equivalent to the
existence of $\{ f_n:k=1,2,\cdots \} $  conjugate functionals  in
the dual Banach space $B^{*}$ of $B$. By $\{f_n\}$ been conjugate
with respect to $ \{e_n\}$ we  means that $f_n(e_m )=\delta_{n m}$
for all $n, m$, where  $\delta_{n m}$ is well know Kronneker's
symbol.
 $\{f_n: n=1, 2,\cdots \}$ is also called a biorthogonal system of
 $E$. It follows that if $E$ is minimal, each $ x \in \overline{\mbox {span}}E$ has a unique formal $E-$expansion $
\sum x_ne_n $ [17], where $ x_n=f_n(x)$.

 In this paper, we particularize $B $ to be the Banach space $
H_{\alpha}$ consisting of all functions $ f(z)$ which are
continuous on the  closed sector $I_{ \alpha}=\{z=re^{i\theta
}:0\leq r \leq 1,|\theta |\leq \alpha \} \ ( 0\leq \alpha <\pi )
$,  analytic in $ \textbf{int}[I_{\alpha}] $. The norm of $f $ is
given by
$$
||f||=\max  \{ |f(z)|: z \in I_{ \alpha }\}.
$$
If the Banach space $H_{\alpha }$ is replaced by
 the Fr$\acute{e}$chlet space $ F_{\alpha } $, which consists of all functions analytic in the sector
 $\textbf{int}[I_{\alpha}]$,
 under the compact topology ( uniform convergence on each compact subset of $\textbf{int}[I_{\alpha}]$ ), Khabibullin [8],
  Rubel [14] and  Malliavin [13] have proved that   $ \textbf{span}[E(\Lambda)]$ is not  dense in
    $ F_{\alpha }  $ if and only if there are $
 b\in (0,\frac{\pi}{\alpha })$ and $M_b $ such that
$$
 \lambda (y)- \lambda (x) \leq b\log y
  - b\log  x +M_b\ \ (y>x \geq 1),
$$
 where
the  characteristic logarithm  $ \lambda(t) $  is defined by ([13]
and [14])
$$
\lambda(t)=\sum_{0<\lambda_n \leq t }\lambda_n^{-1}. \eqno(2)
$$
 Inspired by the method  of  Khabibullin [8],
Anderson [2], Rubel and
 Malliavin, we obtain the following result.

   \vspace{0.3cm}
   \noindent {\bf Theorem 1.}  {\it Let $ \alpha \in
[0,\pi)  $ and  $\Lambda = \{ \lambda _n : n=1,2, \cdots \} $ be a
sequence satisfying $\delta (\Lambda)>0 $ and assume that there
exists a decreasing function $ \varepsilon (x) $ on $[0,\infty ) $
 such that  $ \varepsilon (x)\rightarrow 0$ as $ x\rightarrow \infty $  and such that
$$
 \lambda (y)- \lambda (x) \leq \frac{\alpha}{\pi}\log y
  - \frac{\alpha}{\pi}\log  x + \varepsilon (x),\ \ (y>x \geq 1)
  \eqno(3)
  $$
 Then   $
E(\Lambda) $  is  minimal and $ \overline{\textbf{span}}[
E(\Lambda)] $ does not contain the function $z^{\lambda }$ for  $
\lambda \not \in \Lambda  $, $\mbox{Re}\lambda >0$,     and each
function $ f$ in
 $\overline{ \mbox{span}}( E(\Lambda)) $
  can be extended
 analytically throughout the region $\textbf{int}(I_{\pi})$   with   a  series expansion of the form $(1)$. }

\vspace{0.3cm}
 \noindent
 {\bf Remark 1 } If  $ \lambda(t) $  is  bounded  on $ [1,
\infty )$, then  the  function  $ \varepsilon (x)= \lambda(\infty
)- \lambda(x)$ is decreasing, $ \varepsilon (x)\rightarrow 0$ as $
x\rightarrow \infty $  and   satisfies  (3) with $ \alpha =0$.
 So we have the following Corollary which can be found in [4] ( p.178).

\vspace{0.3cm}
 \noindent
 {\bf Corollary 1. } {\it Let  $\Lambda = \{ \lambda _n : n=1,2, \cdots \} $ be a sequence
satisfying $\delta (\Lambda)>0 $. If  $ \lambda(t) $ is bounded
on $ [1, \infty )$, then     $ E(\Lambda ) \subset H_0$  is
minimal and each function $ f$ in
 $\overline{\textbf{span}}[ E(\Lambda)] $
  can be extended
 analytically throughout the region $\textbf{int}[I_{\pi}]$   with   a  series expansion of the form
 $(1)$. If  $\Lambda = \{ \lambda _n : n=1,2, \cdots \} $ is  a sequence of distinct positive integers and
  $ \lambda(t) $ is bounded  on
$ [1, \infty )$, then each function $ f$ in
 $\overline{ \textbf{span}}[ E(\Lambda)] $
  can be extended
 analytically throughout the open unit disk. }
\noindent
 Therefore,  Theorem 1 is  a generalization   of the Clarkson-Erd\"{o}s-Schwartz  Theorem to a sector.

 \section{Lemmas and Proofs}
In order to prove our main result, we   need  the following
technical
 lemmas. The following Lemma 1 can be seen from [3] and [12].

 \vspace{0.3cm} \noindent
{\bf Lemma 1 ( Fuchs' Lemma ). } {\it If
 $\Lambda  $   is a
  sequence of positive numbers satisfying $\delta( \Lambda )>0$, then
the function
$$
      G(z)=\prod^{\infty}_{n=1} \left (\frac { \lambda_n -z}
      {\lambda_n +z}   \right ) \exp \left (
\frac {2z}{\lambda_n}\right ) \eqno(4)
$$
is  a meromorphic function  and satisfies the following
inequalities:
$$
      |G(z)|\le \exp\{x\lambda (|z|)+Ax\},\ \  z\in \mathbb{C}_{+} ,
     x\geq 0 ,
     \eqno(5)
$$
$$
    |G(z)|\ge \exp\{x\lambda (|z|)-Ax\},\ \  z\in  C(\Lambda ,\delta _0  )
,\eqno(6)
$$
 where  $\  4\delta _0 =\delta(\Lambda) $ and }
$$
  C(\Lambda,\delta_0)=\{z\in \mathbb{C}_+ :
 |z - \lambda_n| \geq \delta_0 , n=1,2,\cdots \}. \eqno(7)
$$

\vspace{0.3cm}

\noindent {\bf Lemma 2. } \hspace{2mm} {\it Let  $ \varepsilon(x)
$ be a positive  decreasing  function on $[0,\infty ) $,  $
\Lambda =\{ \lambda_n:n=1,2, \cdots \} $ and $ \Lambda ^{\prime}
=\{ \lambda_n^{\prime} : n=1,2, \cdots \} $  sequences of positive
numbers
 satisfying  $\delta (\Lambda ) >0 $ and $ \delta (\Lambda ^{\prime}
)>0 $, respectively. If
$$
\lambda (y)-\lambda (x)\leq \lambda ' (y)-\lambda '(x) +
\varepsilon (x), \ y> x\ge 0 , \eqno(8)
$$
 then there exist a constant $A_1$ and  a subsequence $ \Lambda ^* =\{ \lambda_n^*:
n=1,2, \cdots \} $ of the sequence  $ \Lambda ^\prime =\{
\lambda_n^\prime : n=1,2, \cdots \} $   such that
$$
|\lambda (x)+\lambda ^*(x)-
 \lambda ' (x)-A_1|\leq \varepsilon (x)+x^{-1}, \ x>0 . \eqno(9)
$$
}

\vspace{0.3cm}  {\bf Proof of  Lemma 2. }\  Similar to the proof
in [13, p.181-182 ] and [14, p.148-149], we define
$$
   \varphi(x)=\inf \{\lambda '(s)-\lambda (s): s\geq
x  \}.
$$
It follows from (8) that $\varphi (x)\geq \lambda '(x)-\lambda(x)
- \varepsilon (x)$.

 Now $\varphi (x) $ is constant except for possible jumps at the jumps of $\lambda '(x) $.
  Let $a$ be a point of discontinuity of $ \varphi $. Then, the left limit of $ \varphi $ at $ a$ is  $\varphi
  (a-0)=\lambda ' (a-)-\lambda (a-0)$ and the right limit of $ \varphi $ at $ a$ is $\varphi
  (a+0)= \varphi (a)\leq \lambda' (a)-\lambda (a)$. We denote by
  $\Delta \varphi (a)=\varphi
  (a+0)-\varphi  (a-0) $, the jump of $\varphi $ at $a$. Then
 $$
  \Delta
  \varphi (a)\leq \Delta \lambda '(a)- \Delta \lambda (a)\leq
  \Delta \lambda '(a). \eqno(10)
$$
Therefore, there exists a sequence $\Lambda ^* =\{\lambda
^*_n:n=1,2, \cdots \}$ of positive numbers whose counting function
([13] and [14]) is  $ \Lambda ^*(t)=[\Phi(t)]$,
  where $[x]$ denotes the
integral part of $x$,
$$
\Phi (t)=\int_{0}^{t} s  \ d\varphi (s) \ \ \mbox{and }\Lambda
^*(t)=\sum_{\lambda_n^*\leq t }1.
$$
 The characteristic logarithm $ \lambda ^*(t)$ of the
sequence $\Lambda ^* $ is constant except possibly  the jumps of $
\varphi (t)$, and we have $\Delta \lambda^* (a)<a^{-1}+ \Delta
\varphi (a)$. Using (10), we get $\Delta \lambda^* (a)<
a^{-1}+\Delta \lambda '(a)$. Furthermore, $a\Delta\lambda^* (a) $
and $a\Delta \lambda'(a) $ must be integers, so $\Delta \lambda ^*
(a) \leq \Delta \lambda '(a)$ and this means that $ \Lambda ^* $
is a subsequence of $\Lambda ' $.
 Now,
$$
\varphi (x)-\varphi (0)=\int_{0}^{x} s^{-1}  \ d\Phi (s) \ \
\mbox{and }\ \  \lambda ^*(x)=\int_{0}^{x} s^{-1}  \ d[\Phi (s)].
$$
An integration by parts shows that
$$
\lambda^*(x)- \varphi (x)=A_1 +\varepsilon _2(x) ,
$$
where
$$
 \varepsilon _2(x) =\int_{x}^{\infty}(\Phi (s)-[\Phi (s)])\frac{ds}{s^2}- x^{-1}(\Phi
(x)-[\Phi (x)])
$$
and
$$
 A_1=\int_{0}^{\infty} (\Phi (s)-[\Phi (s)])\frac{ds}{s^2}  -\varphi (0).
$$
We now define
$$
\varepsilon_1 (x)=\lambda (x)+\lambda ^* (x)- \lambda '(x)-A_1 ,
 $$
 It is clear that  $ |\varepsilon _2(x)|\leq x^{-1} $, so by the definition of $
   \varphi(x) $,  $ \varepsilon _1(x)\leq \varepsilon _2(x)\leq x^{-1} $. (8)
is simply another way of saying that
$$
\varepsilon_1 (x)\geq \varepsilon_2(x) -\varepsilon (x)\geq
-x^{-1}-\varepsilon (x) .
$$
This proves (9).

\vspace{0.3cm} \noindent {\bf Lemma 3.} {\it Let $ b\geq  0 $ and
let $ \Lambda =\{ \lambda_n:n=1,2, \cdots \} $ be a sequences of
positive real numbers satisfying $\delta (\Lambda )
>0 $.
 If there exists a constant $A_2$ such that
 $$
\lim _{x\rightarrow \infty}|\lambda (x)-b\log^{+} x-A_2|=0,
\eqno(11)
$$
 then the function
$$
 g_0(z)=\frac{G(z)e^{-a_0z}}{\Gamma (\frac{1}{2}+2bz)}, \
  \eqno(12)
$$
 is
meromorphic  and satisfies
$$
\limsup  _{x\rightarrow \infty} x^{-1}\log |g_0(x)|=0, \eqno(13)
$$
$$
   \lim _{x\in C(\Lambda , \delta _0), x\rightarrow \infty} x^{-1} \log|g_0(x)|=0 ,
       \eqno(14)
$$
and
$$
 \lim  _{k\rightarrow \infty} \lambda_k^{-1}\log |g'_0(\lambda
_k)|=0, \eqno(15)
$$
 where  $ \Gamma (z) $ is the Euler Gamma function, $ G(z) $ is
defined by $(4)$ , $a_0=2A_2-2b \log (2b)$ $ ( a_0=2A_2$, if $b=0
)$ and $C(\Lambda ,\delta_0)$ is defined by $ (7)$ }.

\vspace{0.3cm} {\bf Proof of  Lemma 3.}\ The main method of the
proof is based on the use of the following function
   used by Malliavin [11]
       $$
    \psi (s)=2+s\log \left |\frac{ s-1}{s+1}\right | .
$$
 The function  $\psi (s)$ is decreasing on $[0,1)$,  increasing on
$ (1,\infty )$ and there exists $s_0 \in (\frac{5}{6},
\frac{6}{7})$ such that  $ \psi (s_0)=0$. Thus  $\psi  (s) $ is
negative on $ ( s_0, 1  )\bigcup (1,\infty)$. Since $\delta
(\Lambda  )>0$, $ \sum _{n=1}^{\infty}|\lambda _n|^{-2}$
  converges. Thus  $ G(z) $ defined by (4)
 is the quotient of convergent canonical products. As a result,
  the  product (4)
 defines a meromorphic function  in the complex
plane $\mathbb{C}$, which has  zeros  at each point $\lambda_n$.
Writing  $ \log|G(x)| $ as a sum of logarithms, and that sum as a
Stieljes integral, we get
$$
\log |G(x)|=x\int_{0}^{\infty}\psi \left (\frac{t}{x} \right )d
\lambda  (t).
$$
Let
$$
 \  k(x)=\lambda (x)-b\log^{+}
    x-A_2, \ \ \ \varepsilon (x)=\sup \{|k(y)|:y\geq x\}$$
    and
     $$
    \varepsilon _{3}(x)=-\int_{0}^{x}\log \left
|\frac{1-t}{1+t}\right |dt.
$$
Then the function $ \varepsilon _{3}(x) $ is continuous on
$[0,\infty)$, increasing and positive on $(0,\infty) $, convex on
$ [0,1]$ and concave on $[1, \infty)$. Thus $ x^{-1}\varepsilon
_{3}(x) $ is increasing  on
 $ (0,1]$ and decreasing  on $[1, \infty)$, so $  \sup \{ x^{-1}\varepsilon _{3}(x):x>0\}=\varepsilon_3(1)=2\log2<3 $.
 An easy calculation shows that
$$
\int_{0}^{\infty}\psi \left (\frac{t}{x} \right )d \log^{+}t =
\int_{1/x}^{+\infty}\psi (t)\frac{dt}{t}= 2\log x-2+
\varepsilon_3(x^{-1}),
$$
and the Gamma function $ \Gamma (z)$ satisfies
$$
\log \left |\Gamma (\frac{1}{2}+z)\right |=x\log \left
|z+\frac{1}{2}\right | -\left |y\arg (z+\frac{1}{2})\right |-x
+c_1(z),\eqno(16)
$$
where  $ c_1(z)$ satisfies $|c_1(z)| \leq 10 $ for $ x=\mbox{Re}z
\geq 0 $. By the choice of $a_0$, we obtain
$$
x^{-1}\log |g_0(x)|=I_1(x)-2A_2 +\int_{0}^{\infty}\psi \left
(\frac{t}{x} \right )dk(t),\eqno(17)
$$
 where function $ I_1(x)= b\varepsilon_3(x^{-1})-x^{-1}c_1(2bx)$ satisfies
 $$
 \lim_{x\rightarrow \infty}|I_1(x)|=0 .
 $$
Since $k(x)$ has a jump at each  point $\lambda_n$ and   $
\varepsilon (x)\rightarrow 0 $ as $x \rightarrow \infty$,  we can
assume, without loss of generality,  that $1\geq 6\varepsilon
(x)>0$ for $x\ge 0 $( if not, replaced $ \varepsilon (x)$ by $
\min \{\frac{1}{6}, \varepsilon (x)\}$).  Let $a(x)=1+\varepsilon
(\frac{6x}{7})$.
   If we split the range of the integral in (17)
into the ranges $ (0, \frac{x}{a(x)}], [\frac{x}{a(x)}, x a(x)]$
and $ [xa(x),\infty ) $, integration by parts in $ (0,
\frac{x}{a(x)}]$ and $ [xa(x),\infty ) $,  respectively, then  $
x^{-1}\log |g_0(x)| $ can be written in the form
$$
x^{-1}\log |g_0(x)|=\sum _{j=1}^{8}I_j(x),
$$
where
$$ I_2(x)= k\left ( \frac{x}{a(x)}\right )\psi \left (\frac{1}{a(x)}\right ) -k
(xa(x))\psi  (a(x)) ;
$$
$$
I_3(x)=-\int_{0}^{\frac{6x}{7}}k(t)\psi ' \left (\frac{t}{x}\right
)\frac{dt}{x}; \
 \
I_4(x)=-\int_{\frac{7x}{6}}^{\infty}k(t)\psi ' \left
(\frac{t}{x}\right )\frac{dt}{x};
$$
$$
I_5(x)=-\int_{\frac{6x}{7}}^{\frac{x}{a(x)}}k(t)\psi ' \left
(\frac{t}{x}\right )\frac{dt}{x}; \ \ \
I_6(x)=-\int_{xa(x)}^{\frac{7x}{6}}k(t)\psi ' \left
(\frac{t}{x}\right )\frac{dt}{x};
$$
$$
I_7(x)=-b\int_{\frac{x}{a(x)}}^{xa(x)}\psi  \left
(\frac{t}{x}\right )d \log^{+} t \ \ \mbox{and } \
I_8(x)=\int_{\frac{x}{a(x)}}^{xa(x)}\psi  \left (\frac{t}{x}\right
)d \lambda (t).
$$
Next we shall show  that
$$
\lim_{x\rightarrow\infty}|I_j(x)|=0, j=2,3,\cdots,7. \eqno(18)
$$
 Since  $ 1<a(x)=1+\varepsilon \left (\frac{6x}{7}\right )\leq \frac{7}{6}$,
$$
0\leq -\psi \left ( \frac{t}{x}\right )\leq -\log \varepsilon
\left (\frac{6x}{7}\right ) \ \left ( t\in \left [\frac{6x}{7},
\frac{x}{a(x)}\right ]\right  )
$$
  and
$$
0\leq -\psi \left ( \frac{t}{x}\right )\leq -2\log \varepsilon
\left (\frac{6x}{7}\right )\ \ ( t \in [xa(x), \infty )\ ) ,
 $$
we see that
  $$
|I_2(x)|\leq -3 \varepsilon \left (\frac{6x}{7}\right ) \log
\varepsilon \left (\frac{6x}{7}\right );
$$
 $$
 |I_4(x)|\leq -\varepsilon (x)
\int_{\frac{7x}{6}}^{\infty}\psi ' \left (\frac{t}{x}\right
)\frac{dt}{x}=-\varepsilon (x)\psi \left (\frac{7}{6}\right );
$$
$$
|I_5(x)|\leq -\varepsilon \left (\frac{6x}{7}\right
)\int_{\frac{6x}{7}}^{ \frac{x}{a(x)}}\psi ' \left
(\frac{t}{x}\right )\frac{dt}{x}\leq -2 \varepsilon \left
(\frac{6x}{7}\right )\log \varepsilon \left (\frac{6x}{7}\right )
 ;
$$
$$
|I_6(x)|\leq \varepsilon \left (\frac{6x}{7}\right
)\int_{xa(x)}^{\frac{6x}{7})}\psi ' \left (\frac{t}{x}\right
)\frac{dt}{x} \leq - 2\varepsilon \left (\frac{6x}{7}\right )\log
\varepsilon \left (\frac{6x}{7}\right ).
$$
These prove   that (18) hold for $ j=2,4,5,6 $. Also for $x>1$,
 $$ |I_3(x)|\leq
-\varepsilon (0)\int_{0}^{\frac{\sqrt{x}}{2}}\psi ' \left
(\frac{t}{x}\right )\frac{dt}{x} -\varepsilon \left (
\frac{\sqrt{x}}{2}\right
)\int_{\frac{\sqrt{x}}{2}}^{\frac{6x}{7}}\psi ' \left
(\frac{t}{x}\right )\frac{dt}{x}
$$
$$
\leq \frac{\varepsilon (0)}{\sqrt{x}}+2\varepsilon \left
(\frac{\sqrt{x}}{2}\right ),
 $$
so (18) also  holds for $ j=3 .$
 Since $ \frac{1}{a(x)}\geq
\frac{6}{7}>s_0 $, so
$$
0\leq I_7(x)\leq b \varepsilon \left (\frac{6x}{7}\right
)\int_{\frac{x}{a(x)}}^{xa(x)} \left (-\frac{t}{x}\log \left
|\frac{t}{x}-1\right |\right )d \log^{+} t
$$
and $$ 0\leq - I_8(x)\leq \int_{\frac{x}{a(x)}}^{xa(x)} \left
(-\frac{t}{x}\log \left |\frac{t}{x}-1\right |\right )d \lambda
(t).
$$
Therefore
$$
0\leq I_7(x)\leq  2b \varepsilon \left (\frac{6x}{7}\right
)\int_{0}^{\varepsilon (\frac{6x}{7})} (-\log s)ds \leq -4b\left
(\varepsilon \left (\frac{6x}{7}\right )\right )^2\log \varepsilon
\left (\frac{6x}{7}\right ).
$$
Hence (18) holds for $j=7$.  Finally,
$$
0\leq - I_8(x)\leq x^{-1}\sum _{\frac{x}{a(x)}<\lambda_n \leq
xa(x)}\left (\log(3x) - \log |\lambda _n-x |\right).
$$
 let $\Lambda (t)=\sum_{\lambda_n \leq t }1 $ be  the
counting function of  $\Lambda  $ [3], then for $ x\in
C(\Lambda,\delta_0) $, we have $ |\lambda_n-x|\geq
|n-\Lambda(x)|\delta_0 $. Let
$$
n_1(x)=\max \left \{ \Lambda (xa(x))-\Lambda (x),\Lambda
(x)-\Lambda \left( \frac{x}{a(x)}\right ) \right \},
$$
 then, for  $  x\in  C(\Lambda,\delta_0), $
$$
-I_8(x)\leq \frac{2}{x}\left (n_1(x) \log \left(
\frac{3x}{\delta_0 }\right) - \log n_1(x)! \right ).
$$
By  $ e^n n!\geq n^n (n\geq 1),$
$$
-I_8(x)\leq  \frac{2}{x} n_1(x) \left (\log \left ( \frac{3e
x}{\delta_0 }\right) - \log n_1(x) \right ).
$$
 Since the inequalities
 $$
 \Lambda (R)-\Lambda (r)\leq  R (\lambda (R)-\lambda( r))\leq 2R \varepsilon (r)+bR \log \frac{R}{r}  \
 $$
hold for $R>r$,  we see that
  $$
n_1(x)\leq 2xa(x)\varepsilon \left (\frac{6x}{7}\right
)+bxa(x)\log a(x)\leq 4x(1+b)\varepsilon \left (\frac{6x}{7}\right
).
$$
 The function $ t(\log a - \log
t ) $ is increasing on $ (0, ae^{-1}) (a>0)$ and there is $x_0
>1$ such that $9\delta_0 (1+b)\varepsilon
\left (\frac{6x}{7}\right )\leq 3$ , we see that
 $$
 -I_8(x)\leq -18(1+b) \varepsilon
\left (\frac{6x}{7}\right ) \log \left ( \delta_0(1+b)\varepsilon
\left (\frac{6x}{7}\right )\right ) , \ \ x\ge x_0 .
$$
 These prove that (13) and (14) hold. Similarly, (15) can also be
 proved.
This completes the proof of Lemma 3. \vspace{0.3cm}

\section{Proof of Theorem}
 {\bf Proof . }\ We can assume that $ \alpha>0 $ in the proof of
Theorem.
 It is  a
consequence of the Hahn-Banach theorem [15] that
$\overline{\textbf{span}}[E(\Lambda )]\neq H_{\alpha}$  if and
only if there
 exists a bounded linear functional $T$ on $ H_{\alpha}$  with  $ ||T||=1 $ which vanishes
 on all of  $E(\Lambda
)$. Since every bounded linear functional on $ H_{\alpha}$  is
given by integration with respect a complex Borel measure on $
I_{\alpha }$.  So   we shall construct a bounded linear functional
$ T$ on $H_{\alpha}$ such that
\begin{equation*}
T(\zeta^{z})= g(z)=\frac{z^2G(z)e^{-Az}}{\Gamma
(\frac{1}{2}+\frac{2}{\pi}\alpha z)(1+z)^{4}},
\end{equation*}
 where  $ \Gamma (z) $ is the Euler Gamma function, $ G(z) $ is
defined by (4) and  $A $ is a sufficient large positive constant.
The function $ g(z) $ is analytic in the right half plane $
\mathbb{C}_{+}$. Moreover, since $ G(z)G(-z)\equiv 1$ and $\Gamma
(z)\Gamma (1-z)\sin (\pi z) \equiv \pi$,  it follows from Lemma 3,
(16) and Cauchy's formula for $g'(z) $ and $g''(z)$  that
$$
   |g(z)|+|g'(z)|+ |g''(z)| \leq \frac {Ae^{ \alpha |y| }}{1+|z|^2} \ \  (x\geq 0) \eqno(19)
$$
 holds for a sufficient large  positive  constant $ A$.
   Fix $ z $ so that $ x>0, y>0$, and consider the Cauchy's  formula for
   $ g(z)e^{\alpha z i} $, where the path of integration consists of
   the  quadrant circle with center at $0 $, radius $ R>1+|z|$ from
   $R$ to $iR$, followed by the interval from $iR $ to $0$ and by
   the  interval from $0$ to $ R$. The integral over the quadrant circle
tends to $ 0$ as $ R\rightarrow \infty $, so we are left with
$$
g(z)e^{i\alpha z} =\frac {1}{2\pi i} \int ^{+\infty}_{0}
\frac{g(t)e^{i\alpha t}}{t-z}dt - \frac {i}{2\pi i} \int
^{+\infty}_{0} \frac{g(it)e^{i\alpha it}}{it-z}dt  \eqno(20)
$$
and similarly,  fix $ z $ so that $ x>0, y>0$, and consider the
Cauchy formula for
   $ g(z)e^{-i\alpha z } $, where the path of integration consists of
   the lower  quadrant circle with center at $0 $, radius $ R>1+|z|$ from
   $-iR$ to $R$, followed by the interval from $R $ to $0$ and by
   the  interval from $0$ to $ -iR$. The integral over the lower quadrant circle
tends to $ 0$ as $ R\rightarrow \infty $, so we are left with
$$
0 =\frac {-1}{2\pi i} \int ^{+\infty}_{0} \frac{g(t)e^{-i\alpha
t}}{t-z}dt + \frac {i}{2\pi i} \int ^{-\infty}_{0}
\frac{g(it)e^{-i\alpha it}}{it-z}dt. \eqno(21)
$$
Using
$$
\frac{1}{z-it}=\int_0^1 s^{z-it-1}ds \ \ \mbox{and} \ \
\int_{-\alpha}^{\alpha }e^{i\theta (t-z)}id\theta
=\frac{e^{i\alpha (it-z)}-e^{-i\alpha (it-z)}}{t-z}
$$
 (20) multiplied by $e^{-\alpha z i}$ plus (21) multiplied by $e^{\alpha z i}$, we   obtain, for $ z=x+iy, x>0,y>0 $,
$$
g(z) =\frac {1}{2\pi } \int ^{+\infty}_{0}
g(t)\int_{-\alpha}^{\alpha }e^{i\theta (t-z)}d\theta dt
$$
\begin{equation*}
- \frac{1}{2\pi } \int _{-\infty}^{0} g(it)\int_0^1
(se^{i\alpha})^{ (z-it)}\frac{ds}{s}dt -\frac{1}{2\pi } \int
^{+\infty}_{0} g(it)\int_0^1 (se^{-i\alpha})^{
(z-it)}\frac{ds}{s}dt.\eqno(22)
\end{equation*}
Similarly, (22) also holds for $x>0, y<0$. The interchange in
order of integration in (22) are legitimate: in the integrants in
(22) are replaced by their  absolute values, some finite integral
results. Hence (22) can be rewritten in the form
\begin{equation*}\label{25}
g(z) =\frac {1}{2\pi } \int_{-\alpha}^{\alpha }e^{i\theta
z}h_0(e^{i\theta}) d\theta
\end{equation*}
\begin{equation*}
 +\frac{1}{2\pi i} \int_0^1 ((se^{-i\alpha})^{
z}h_1(se^{-i\alpha})-(se^{i\alpha})^{z}h_{-1}(se^{i\alpha}))\frac{ds}{s},
\end{equation*}
where
$$
   h_l(\zeta)=\int _{L_l} g(z)
   \zeta^{-z}dz \eqno(23)
$$
 and  $L_l=\{t\exp \{ \frac{\pi}{2}li \}: t\geq 0 \}\ \  (l\in
   \{-1,0,1\}) $ are half-lines.
   By (21),
$h_0(\zeta ) $ is  analytic in the region  $ D_0=\{\zeta : |\zeta|
>1, |\arg \zeta |<\pi \}$ and continuous in the set $\overline{ D_0}=\{\zeta : |\zeta| \geq
1 , |\arg \zeta |<\pi \}$,  each function  $h_l(\zeta ) (l=\pm 1)
$ is  analytic in the sector $ D_l=\{\zeta : \alpha <-l\arg \zeta
<\pi \}$ and continuous in the closure $\overline{D_l}= \{\zeta :
\alpha \leq -l\arg \zeta \leq \pi \}$ of $D_l$.
 By Cauchy's
formula, $ h_{0}(\zeta ) $ can be continued analytically to a
bounded analytic function in the region $D_{-1}\bigcup D_0\bigcup
D_1$ $ = \{\zeta=\rho e^{i\phi }: \zeta \notin I_{\alpha }, |\phi
|<\pi \}$, i.e., $ h_0(\rho e^{i\phi} )=h_{l}(\rho e^{i\phi })$
for $ \rho
>1, \alpha <- l \phi <\pi, l=\pm 1 $.  By (21), $h_0(\zeta) $ is
bounded in the circular arc $ \{\zeta: |\zeta|=1, |\arg \zeta
|<\alpha \} $. Integrations by parts twice in (23),
\begin{equation*}\label{24}
   h_l(se^{-il \alpha })=(\log s -il \alpha )^{-2}\int _{L_l} g''(z)
   (se^{-il\alpha })^{-z}dz , \ \ l=\pm 1 .
 \end{equation*}
By (19),
\begin{equation*}\label{25}
\int_0^1
(|h_{-1}(se^{i\alpha})|+|h_1(se^{-i\alpha})|)\frac{ds}{s}< \infty.
\end{equation*}
Therefore, the linear functional
\begin{equation*}\label{26}
T(\varphi ) =\frac {1}{2\pi } \int_{-\alpha}^{\alpha }\varphi
(e^{i\theta })h_0(e^{i\theta}) d\theta
\end{equation*}
\begin{equation*}
  +\frac{1}{2\pi i} \int_0^1 (\varphi
(se^{-i\alpha})h_1(se^{-i\alpha })-\varphi
(se^{i\alpha})h_{-1}(se^{i\alpha} ))\frac{ds}{s}
\end{equation*}
is a bounded linear functional on $H_{\alpha }$ and
 satisfies  $ T(\zeta^{\lambda })=g(\lambda ) $ for
$\lambda \in \mathbb{C}_{+} $.
 By the Riesz representation theorem,  $ z^{\lambda} \not\in \overline{\mbox{span}}E(\Lambda
)$ for $\lambda \not \in \Lambda $ and $\mbox{Re} \lambda >0 $.
Similarly, replacing $g(z)$ by $ (z-\lambda )^{-1}g(z)$ for $
\lambda \in \Lambda $, we can also prove that no element of
$E(\Lambda ) $ belongs to the closure of the vector subspace
generated by the other elements of $E(\Lambda )$. Therefore,
$E(\Lambda ) $ is minimal, and $\overline{\mbox{span}} M(\Lambda
)\neq H_{\alpha } $.

 Next,  define, for $0 <b=\frac{\alpha}{\pi}< \infty $, the arithmetic progression [13] $\
\Lambda_b \ $  by
\begin{equation*}
\Lambda_b=\left \{ \frac{n}{b}:n=1,2,\cdots \right \}
\end{equation*} and observe that the counting function $ \Lambda_b(t)=\sum _{n\leq b}1 $
of $ \Lambda_b $ satisfies $ \Lambda_b(t)=[bt]=bt+O(1),$ and the
characteristic logarithm $ \lambda _b(t)$ of the sequence
$\Lambda_b $ satisfies
\begin{equation*}
  \lambda_b(t)=b\log t+b\log b + b\gamma
+O(t^{-1}),
\end{equation*}
$ \mbox {as} \ t\rightarrow \infty $, where $\gamma $ is a Euler
constant. So by
 Lemma 2, there exist a constant $A_1$ and  a subsequence $ \Lambda ^* =\{ \lambda_n^*:
n=1,2, \cdots \} $ of $ \Lambda_b $ such that (9) holds. If $
\Lambda $ and $  \Lambda ^*  $ have common elements or $\delta (
\Lambda \cup  \Lambda ^*)=0$, we adjust $ \Lambda ^* $ as follows:
let $ 4h_1=min \{ \delta (\Lambda ), b  \} $ and $n_{k}\in
\mathbb{N}$ such that $\lambda_{n_{k}}\leq\lambda_{k}^{*}\leq
\lambda _{1+n_{k}}$ and let
\[
\lambda_{k}^{**}=\left\{
\begin{array}{ll}
\lambda_{k}^{*},&\quad \mbox{if} \
\lambda_{n_k}+h_1\leq\lambda_{k}^{*}<\lambda_{1+n_k}-h_1;\\
\lambda_{k}^{*}+h_1,&\quad \mbox{if} \ \lambda_{n_k}\leq\lambda_{k}^{*}<\lambda_{n_k}+h_1;\\
\lambda_{k}^{*}-h_1,&\quad \mbox{if } \
\lambda_{1+n_k}-h_1<\lambda_{k}^{*}<\lambda_{1+n_k},
\end{array}
\right.
\]
and let
$$
A_3=\sum_{k=1}^{+\infty}\left
(\frac{1}{\lambda_{k}^{*}}-\frac{1}{\lambda^{**}_{k}}\right ),
$$
then the set $\Lambda^{**} =\{ \lambda_n^{**}:n=1,2,\cdots \}$ and
 the set  $ \Lambda  $ are  disjoint  and $\delta ( \Lambda
\cup \Lambda ^{**})\geq  h_1>0 .$  For $x\geq 2h_1+1 $, we have
the following inequalities:
$$
|\lambda^{*}(x)-\lambda^{**}(x)-A_3|\leq\frac 1 x
+\sum_{\lambda^{*}_{k}\geq x}^{} \frac{h_1}{\lambda
_k\lambda_k^*}\leq\frac {1}{x}+\frac{1}{x-h_1}\leq\frac{3}{x}
$$
and
$$
|\lambda (x)+\lambda^{**}(x)-\frac{\alpha }{\pi}\log x -A_1-
A_3|\leq \frac {13}{x}+\varepsilon (x).
$$
  Suppose that $ f$ is in
  $ \overline{\textbf{span}}[E(\Lambda)]$, since
  $\overline{\textbf{span}}[E(\Lambda)]
   \subset \overline{\textbf{span}}[E(\Lambda \cup \Lambda^{**})],$
then from the uniqueness of $E(\Lambda \cup
\Lambda^{**})-$expansion $\sum b_nz^{\mu_n}$  of $f$,
 where  $\ \Lambda \cup \Lambda^{**}=\{ \mu _n:n=1,2,\cdots \}$, we see that those
 coefficients $b_n $ associated with members $ \mu _n\in \Lambda^{**} $ distinct from
 all $\lambda_n $ are equal to zero. Thus the $E(\Lambda \cup \Lambda^{**})-$expansion reduces to
 $E(\Lambda )-$expansion whenever $f \in \overline{\mbox{span}}E(\Lambda )$. Therefore,
  we can assume, without loss of generality,  that there exists a constant $ A_2$ such that (11) holds
  with $ b=\frac{\alpha}{\pi }$. Therefore, the function $g_0(z)$ defined by (12) satisfies (13),(14) and (15).
Let
 $$
 \psi_k(z)=\frac{z^2g_0(z)}{(1+z)^{4}(z-\lambda_k)},\ \
 \psi_k(\lambda_k)=\frac{\lambda_k^2g'_0(\lambda_k)}{(1+\lambda_k)^{4}},
 $$
and
\begin{equation*}
 h_{k,l}(\zeta )=\int _{L_l} \psi_k(z)
   \zeta^{-z}dz ,
 \end{equation*}
where  $L_l=\{t\exp \{ \frac{\pi}{2}li \}: t\geq 0 \} (l\in
   \{-1,0,1\}) $ are half-lines.  As has been shown in (19), then there exists a positive constant $A_4$ such that
$$
  |\psi_k(iy)|+|\psi_k'(iy)|+|\psi_k''(iy)| \leq \frac {A_4e^{ \alpha |y|
  }}{1+|y|^2} \eqno(24)
$$
and
$$
  \limsup  _{x\rightarrow \infty} x^{-1}\log |\psi_k(x)|=0 \eqno(25)
$$
hold for each $k$. By (15),
$$
  \limsup  _{k\rightarrow \infty} \lambda_{k}^{-1}\log |\psi_k(\lambda _k)|=0. \eqno(26)
$$
By (24) and (25), $h_{k,0}(\zeta ) $ is analytic in the region $
D_0=\{\zeta : |\zeta|
>1, |\arg \zeta |<\pi \}$ ,   $h_{k,l}(\zeta )
(l=\pm 1) $ is  analytic in the sector $ D_l=\{\zeta : \alpha
<-l\arg \zeta <\pi \}$ and continuous in the closure
$\overline{D_l}= \{\zeta : \alpha \leq -l\arg \zeta \leq \pi \}$
of $D_l$.
 By Cauchy's
formula, $ h_{k,0}(\zeta ) $ can be continued analytically to an
analytic function in the region $D_{-1}\bigcup D_0\bigcup D_1$ $ =
\{\zeta=\rho e^{i\phi }: \zeta \notin I_{\alpha }, |\phi |<\pi
\}$, i.e., $ h_{k,0}(\rho e^{i\phi} )=h_{k,l}(\rho e^{i\phi })$
for $ \rho
>1, \alpha <- l \phi <\pi, l=\pm 1 $.  By (25), $h_0(e^{-\delta}\zeta) $ is
bounded in the circular arc $ \{\zeta: |\zeta|=1, |\arg \zeta
|<\alpha \} $ for each $\delta >0$.
   the  linear functionals
\begin{equation*}\label{26}
T_{k,\delta }(\varphi ) =\frac {1}{2\pi } \int_{-\alpha}^{\alpha
}\varphi (e^{i\theta })h_{k,0}(e^{-\delta}e^{i\theta}) d\theta
\end{equation*}
\begin{equation*}
+\frac{1}{2\pi i} \int_0^1 (\varphi
(se^{-i\alpha})h_{k,1}(e^{-\delta}se^{-i\alpha })- \varphi
(se^{i\alpha})h_{k,-1}(e^{-\delta}se^{i\alpha}  )\frac{ds}{s}\ \
(\varphi \in H_{\alpha })
\end{equation*}
are  bounded linear functionals in $ H_{\alpha }$ and  satisfy  $
T_{k,\delta }(\zeta ^{\lambda })=\psi_k(\lambda )e^{-\delta
\lambda} $ for $\lambda \in \mathbb{C}_+$ and
$A(\delta)=\sup\{||T_{k,\delta}||:k=0,1,2, \cdots\}<\infty $.
Therefore, $ \{e^{\delta \lambda_k
}(\psi_k(\lambda_k))^{-1}T_{k,\delta}: k=1,2,\cdots \} $ is a
biorthogonal system of $E(\Lambda) $.
 If $f$ belongs to  $ \overline{\textbf{span}}[E(\Lambda)]$,
then there exists a sequence of $\Lambda -$polynomials
$$
 P_l(z)=\sum_{n=1}^{l}a_{n,l}z^{\lambda_n } \in
\mbox{span}E(\Lambda )
$$
 such that
$$
||f-P_l|| \longrightarrow 0\ \ {\rm as}\ \ l\longrightarrow
\infty.
$$
 Let
  $$
   \sum_{k=1}^{\infty}a_{k}z^{\lambda_k}
\eqno(27)
$$
be the $E(\Lambda )-$expansion  of $f$.   The  biorthogonality of
the system
 $$
 \{e^{\delta \lambda_k
}(\psi_k(\lambda_k))^{-1}T_{k,\delta}:k=1,2,\cdots \} $$
 implies
that
 $ a_k =e^{\delta \lambda _k}(\psi_{k}(\lambda_k
))^{-1}{T}_{k,\delta}(f) $ and $
 a_{k,l}=e^{\delta \lambda _k}(\psi_{k}(\lambda_k ))^{-1}{T}_{k,\delta}(P_l) $. Therefore
$$
|a_k-a_{k,l}| \leq ||f-P_l||A(\delta )e^{\delta \lambda
_k}|\psi_{k}(\lambda_k )|^{-1} \ \ ( k=1,2,\cdots)
$$
(thus  that the sequences $\{ a_{k,l}\}$ are independent of
$\delta $ implies that the sequence $\{a_{k}\}$ is  also
independent of $\delta $ ) and
  $$  |a_{k}| \leq
 A(\delta )||f||e^{\delta \lambda _k}|\psi_{k}(\lambda_k )|^{-1}, \ \  k=1,2,\cdots . $$
By (26),   the   series in (27) converges to
  an analytic function $F(z)$ uniformly on compacts of  $\{z: |z|<1, |\arg z|<\pi \} $.
  we obtain that, for $ z\in \textbf{int}[I_{\alpha}]$, there is $\delta>0$ such that $|z|<e^{-\delta}
  $, so
\begin{align*}
    |f(z)-F(z)| & \leq |f(z)-P_l(z)|+|P_l(z)-F(z)| \\
& \leq ||f-P_l||+\sum_{n=1}^{l} |a_{nl}-a_n||z|^{ \lambda _n
}+\sum_{n=l+1}^{\infty}|a_n||z|^{\lambda_n  }.
\end{align*}
Letting
 $l\rightarrow \infty $, we obtain that $ f(z)=F(z) $ for
$ z\in $int$I_{\alpha}$. This completes the proof of Theorem.


\begin{thebibliography}{10}
\bibitem{1}J. M. Almira, M\"{u}ntz type theorems, \emph{Surveys in Approximation Theory} {\bf 3}(2007), 152-194.

\bibitem{2}J.M. Anderson, M\"untz-Sz\'asz type
approximation and the angular growth of lacunary integral
functions, \emph{ Trans. Amer. Math. Soc.}, {\bf 169}(1972),
237--248.


\bibitem{3}R. P. Boas, Jr.,   Entire Functions,  \emph{ Academic Press, New
York }, (1954).


 \bibitem{4} P. B. Borwein and T. Erd\'{e}lyi,  Polynomials and
 Polynomial Inequalities, \emph{ Springer-Verlag, New York}, (1995)

\bibitem{5}J. A. Clarkson, P.Erd\"{o}s,  Approximation by
polynomials, \emph{Duke Math.J.} {\bf 10}(1943), 5--11.

\bibitem{6} T. Erd\'{e}lyi, The full
Clarkson-Erd\"{o}s-Schwartz  theorem on  the closure of non-dense
M\"{u}ntz spaces, \emph{ Studia Math. } {\bf 155}(2003), 145-152.


\bibitem{7} V. I. Gurariy, Geometry of  M\"{u}ntz spaces and
related questions, Lecture Notes in Mathematics  Vol {\bf 1870},
\emph{Springer,} (2005).

\bibitem{8}B.N. Khabibullin, On the growth of entire functions
of exponential tye along the imaginary axis, \emph{ Soviet Math.
Dokl.}, {\bf 38}(2)(1989), 276--278 .

\bibitem{9}B.N. Khabibullin, \ \ On the growth of entire functions
of exponential type along the imaginary axis, \emph{ Math. USSR
Sbornik }, {\bf 67}(1)(1990), 149--163.

\bibitem{10} P. Koosis,  The Logarithmic integral II,
\emph{Cambridge University Press}, (1992)

\bibitem{11} P. Malliavin, Sur la croissance radiale d'une
 fonction m\'{e}romorphe, \emph{Illinois J. of Math.},
{\bf 1}(1957), 259--296.

\bibitem{12}P. Malliavin, Sur quelques proc\'ed\'es
d'extrapolation, \emph{Acta Math.}, {\bf 83}(1955), 179--255.

\bibitem{13}P. Malliavin and L.A. Rubel, On small entire
functions of exponential type with given zeros, \emph{ Bull. Soc.
Math. France,} {\bf 89}(1961), 175--206.



\bibitem{14}L.A. Rubel, Entire and meromorphic funcions,
\emph{Springer-Verlag, New York, Inc.}, (1996).


\bibitem{15}W. Rudin, Real And Complex Analysis,  \emph{
McGraw-Hill, New York}, (1974).

\bibitem{16} L. Schwartz, \'{E}tude des sommes d'exponentielles,
\emph{ Hermann, Paris }, (1959).


 \bibitem{17}
I. Singer,\ \
 Bases in Banach Spaces I. \emph{ Springer-Verlag, New
York}, (1970).

\bibitem{18} O. Sz\'{a}sz, \"Uber die Approximation steliger
Funktionen durch lineare Aggregate von Potenzen, \emph{Math.
Ann.}, {\bf 77}(1916), 482--496.






\end{thebibliography}
\end{document}